\documentclass[12pt,a4paper,twoside]{article}
\usepackage[english]{babel}
\usepackage{amssymb}
\usepackage{inputenc}
\usepackage{amssymb}
\begin{document}

\author{S. Albeverio $^{1},$ Sh. A. Ayupov $^{2},$ A. A. Zaitov $^3,$ J. E. Ruziev $^4$}

\title{\bf Algebras of unbounded operators over the ring of measurable
functions and their derivations and automorphisms}

\maketitle

\begin{abstract}

In the present paper derivations and $*$-automorphisms of algebras
of unbounded operators over the ring of measurable functions are
investigated and it is shown that all $L^0$-linear derivations and
$L^{0}$-linear $*$-automorphisms are inner. Moreover, it is proved
that each $L^0$-linear automorphism of the algebra of all linear
operators on a $bo$-dense submodule of a Kaplansky-Hilbert module
over the ring of measurable functions is spatial.

\end{abstract}

\medskip
$^1$ Institut f\"{u}r Angewandte Mathematik, Universit\"{a}t Bonn,
Wegelerstr. 6, D-53115 Bonn (Germany); SFB 611, BiBoS; IZKH; HIM;
CERFIM (Locarno); Acc. Arch. (USI), \emph{albeverio@uni-bonn.de}

$^2$ Institute of Mathematics and information technologies,
Uzbekistan Academy of Science, F. Hodjaev str. 29, 100125, Tashkent
(Uzbekistan), e-mail: \emph{sh\_ayupov@mail.ru,
e\_ayupov@hotmail.com, mathinst@uzsci.net}

 $^{3}$ Institute of
Mathematics and information technologies, Uzbekistan Academy of
Science, F. Hodjaev str. 29, 100125, Tashkent (Uzbekistan), e-mail:
\emph{adilbek$_{-}$zaitov@mail.ru}

$^{4}$ Institute of Mathematics and information technologies,
Uzbekistan Academy of Science, F. Hodjaev str. 29, 100125, Tashkent
(Uzbekistan), e-mail: \emph{jaloler@mail.ru}



\newpage
\large

\section*{\center 0. Introduction}

The theory of derivations and automorphisms of operator algebras is
an important branch of the theory of operator algebras and
mathematical physics. The present paper is devoted to the study of
derivations and automorphisms of the algebras of unbounded operators
over the ring of measurable functions. Derivations on the algebras
of bounded operators are rather well-investigated [1]. A certain
method of investigation of derivations was suggested in [3], where
it was proved that any derivation of a standard algebra of bounded
operators on a normed space is inner and any automorphism of such
algebra is spatial.

A survey of results and open problems in the theory of derivations
on unbounded operators algebras were given in [2]. Later the
existence of non-inner derivations on the algebra $L(M)$ of
measurable operators affiliated with an abelian von Neumann algebra
$M$ was established in [4]. Recently it was proved [5] that in the
algebra of (equivalence classes of) measurable complex functions on
a locally separable measure space there exist non trivial
derivations and non-extendable automorphisms which are not
identical.

Derivations and automorphisms of special classes of unbounded
operator algebras (so-called $O^\ast$-algebras) were considered in
[6], in particular it was proved that all derivations and all
$*$-automorphisms of the maximal $O^*$-algebra
$\mathcal{L}^{+}(\mathcal{D})$ are inner and every automorphism of
the algebra $\mathcal{L}(\mathcal{D})$ is spatial. In the present
paper we study derivations and automorphisms of standard algebras of
unbounded $L^0$-linear operators and obtain $L^0$-valued versions of
the above results from [6].

It should be noted that $L^0$-valued analogues of some classic
results become very useful in solving problems of classical operator
algebras. For example, in [7] the theory of Kaplansky-Hilbert
modules over $L^0$ has been applied for the investigation of
derivations on algebras of $\tau$-measurable operators affiliated
with a type $I$ von Neumann algebra and faithful normal semi-finite
trace $\tau.$

The Section 1 contains preliminaries from the theory of
Kaplansky-Hilbert modules over $L^{0}$. In Section 2 we develop the
theory of unbounded $L^{0}$-linear operators on Kaplansky-Hilbert
modules over $L^0$ and introduce and study notions such as
$O$-modules, $O^*$-modules, $O$-algebras, $O^*$-algebras for the
$L^0$-valued case. Further we show that every $L^0$-linear
derivation of the algebra $\mathcal{L}^{+}(\mathcal{D})$ is inner
and each automorphism of the algebra $\mathcal{L}(\mathcal{D})$ is
spatial. We also consider $*$-isomorphisms of $O^*$-algebras over
the ring of measurable functions and prove that every $L^0$-linear
$*$-isomorphism between $O^*$-algebras is spatial and each
$L^0$-linear $*$-automorphism of the algebra
$\mathcal{L}^{+}(\mathcal{D})$ is inner.

\section*{\center 1. Kaplansky-Hilbert modules over the ring of measurable functions}

Let $(\Omega,\Sigma,\mu)$ be a space with a complete finite measure,
and let $L^0=L^{0}(\Omega)$ be the algebra of all measurable
complex-valued functions on $(\Omega,\Sigma,\mu)$ (functions equal
almost everywhere are identified).

Consider a vector space $X$ over the field $\mathbb{C}$ of complex
numbers. A map $\|\cdot\|:X\longrightarrow L^0$ is called an
$L^0$-valued norm on $X$, if for any
$\varphi,\psi\in X,$ $\lambda\in \mathbb{C}$ the following conditions are fulfilled:\\
$1) \|\varphi\|\geq 0;\\
2) \|\varphi\|=0\Longleftrightarrow \varphi=0;\\
3) \|\lambda \varphi\|=|\lambda|\|\varphi\|;\\
4) \|\varphi+\psi\|\leq\|\varphi\|+\|\psi\|.$

The pair $(X,\|\cdot\|)$ is said to be a \emph{lattice-normed} space
(shortly, LNS) over $L^0$. An LNS $X$ is called $d$-decomposable, if
for any $\varphi\in X$ and for each decomposition
$\|\varphi\|=e_{1}+e_{2}$ into the sum of disjoint elements there
exist $\varphi_{1}, \varphi_{2}\in X$ such that
$\varphi=\varphi_{1}+\varphi_{2}$ and $\|\varphi_1\|=e_{1},
\|\varphi_2\|=e_{2}$. A $d$-decomposable norm is also called a
\emph{Kantorovich norm}. A net $(\varphi_\alpha)_{\alpha\in A}$ of
element from $X$ is called $(bo)$-{\it convergent} to $\varphi\in
X$, if the net $(\|\varphi_{\alpha}-\varphi\|)_{\alpha\in A}$
$(o)$-converges to zero in $L^0$ (recall that $(o)$-convergence of a
net from $L^0$ is equivalent to its convergent almost everywhere). A
{\it Banach-Kantorovich space} (further, BKS) over $L^0$ is a
$(bo)$-complete $d$-decomposable LNS over $L^0$.

Any BKS $X$ over $L^0$ is a module over $L^0$, i. e. for any
$\lambda\in L^0$ and $\varphi\in X$ the element $\lambda\varphi\in
X$ is determined and $\|\lambda \varphi\|=|\lambda|\|\varphi\|$ (see
[8, 9]).

A module $E$ over $L^0$ is said to be \emph{finite-generated,} if
there exist $\varphi_{1},\varphi_{2},...,\varphi_{n}$ in $E$ such
that every $\varphi\in E$ can be decomposed as
$\varphi=\alpha_{1}\varphi_{1}+...+\alpha_{n}\varphi_{n}$ where
$\alpha_{i}\in L^0, i=\overline{1,n}$. The elements
$\varphi_{1},\varphi_{2},...,\varphi_{n}$ are called
\textit{generators} of the module $E.$ A minimal number of
generators of a finite-generated module $E$ is denoted by $d(E).$ A
module $E$ over $L^0$ is called \emph{$\sigma$-finite-generated,} if
there exists a partition $(\pi_{n})_{n\in\mathbb{N}}$ of the unit in
$\nabla$ ($\nabla$ is the Boolean algebra of all idempotents in
$L^0$) such that each $\pi_{n}E$ is finite-generated. A
finite-generated module $E$ over $L^0$ is called \textit{homogeneous
of type} $n,$ if $n=d(\pi E)$ for every nonzero $\pi\in\nabla$.

Elements $\varphi_{1},\varphi_{2},...,\varphi_{n}\in E$ are called
\emph{$\nabla$-linear independent}, if for every $\pi\in\nabla$ and
any $\alpha_{1},\alpha_{2},...,\alpha_{n}\in L^0$ the equality
$\pi\sum\limits_{k=1}^{n}\alpha_{k}\varphi_{k}=0$ implies
$\pi\alpha_{1}=\pi\alpha_{2}=...\pi\alpha_{n}=0$ (see [7]).

If $E$ is module over $L^0$ which is a homogeneous of type $n$ then
there exists a basis $\{\varphi_{1},\varphi_{2},...,\varphi_{n}\}$
in $E,$ consisting of $\nabla$-linear independent elements, i. e.
each element $\varphi\in E$ can be uniquely represented in the form
$\varphi=\alpha_{1}\varphi_{1}+...+\alpha_{n}\varphi_{n}$,
$\alpha_{i}\in L^0, i=\overline{1,n}$ (see [10], Proposition 6).

Let $X$ and $Y$ be BKS over $L^0.$ An operator $a:X\rightarrow Y$ is
$L^0$-{\it linear} if $a(\alpha \varphi +\beta \psi )=\alpha
a(\varphi)+\beta a(\psi)$ for all $\alpha,\beta \in L^0 ,\varphi ,
\psi \in X.$ The set of all $L^0$-{\it linear} operators is denoted
by $\mathcal{L}(X, Y).$ An operator $a\in \mathcal{L}(X, Y)$ is
called $L^0$-\textit{bounded}, if there exists $c\in L^0$ such that
$\|a(\varphi)\|\leq c\|\varphi\|$ for all $\varphi\in X.$ For an
$L^0$-bounded operator $a$ we put
$\|a\|=\sup\{\|a(\varphi)\|:\|\varphi\|\leq \textbf{1} \}.$ An
$L^0$-linear operator $a:X\rightarrow Y$ is said to be
\textit{finite-generated} (respectively,
$\sigma$-\textit{finite-generated}, \textit{homogeneous of type}
$n$), if $a(X)=\{a(\varphi):\varphi\in X\}$ is a finite-generated
(respectively, $\sigma$-finite-generated, homogeneous of type $n$)
submodule in $Y.$

It is clear that each $ L^0$-linear $\sigma$-finite-generated
operator $a:X\rightarrow Y$ can be represented as
$a=\sum\limits_{n=1}^{\infty}\pi_{n}a_{n},$ where
$(\pi_{n})_{n\in\mathbb{N}}$ is a partition of the unit $\nabla,$
and $a_n$ are homogeneous operators of finite type. Moreover if $a$
is a finite-generated operator then $(\pi_{n})$ is a finite
partition of unit.

Let $a:X\rightarrow Y$ be a homogeneous of type $n$ $ L^0$-linear
operator and let $\{\psi_1,..., \psi_n\}$ be a basis in $a(X).$
Denote by $X^{\ast}$ the space of all $ L^0$-bounded $ L^0$-linear
functionals from $X$ into $ L^0.$ Then there exists a system
$\{f_1,..., f_n\}\subset Y^{\ast}$ such that $f_i(\psi_j)=\delta_{i
j}\textbf{1},$ where $\delta_{i j}$ is Kroenecker symbol (see [10],
Proposition 2). We define $g_i\in X^{\ast}, i=\overline{1, n}$ as
follows
$$g_i(\varphi)=f_i(a(\varphi)),\quad \varphi\in X.$$
It is clear that
$$a(\varphi)=\sum\limits_{k=1}^{n}g_{k}(\varphi)\psi_{k},\quad \varphi\in X.$$
This formula gives the general form of  $L^0$-bounded $ L^0$-linear
operators from $X$ into $Y$ which are homogeneous of type $n
(n\in\mathbb{N})$.

If $X$ and $Y$ coincide then $\mathcal{L}(X)$ is used for
$\mathcal{L}(X, X).$

An algebra $\mathcal{U}\subset\mathcal{L}(X)$ over $ L^0$ is said to
be \textit{standard} if $\mathcal{F}(X)\subset\mathcal{U}$, where
$\mathcal{F}(X)$ is the algebra of all finite-generated $
L^0$-linear operators from $\mathcal{L}(X)$. The following algebras
over $L^0$ are examples of standard algebras: the algebra
$\mathcal{F}(X)$; the algebra $\mathcal{F}_{\sigma}(X)$ of all
$\sigma$-finite-generated $L^0$-linear operators from
$\mathcal{L}(X)$; the algebra $\mathcal{K}(X)$ of all $L^0$-linear
cyclically compact operators from $\mathcal{L}(X)$; the whole
algebra $\mathcal{L}(X).$

Let $\mathcal{A}$ be a module over $ L^0$. A map $\langle
\cdot,\cdot\rangle:\mathcal{A}\times \mathcal{A}\rightarrow L^0$ is
called an $ L^0$-\textit{valued inner product}, if for all
$\varphi,\psi,\eta\in \mathcal{A},\,\lambda\in L^0$ the following conditions are fulfilled:\\
$1) \langle \varphi,\varphi\rangle\geq0;\\
2) \langle \varphi,\varphi\rangle=0\Leftrightarrow \varphi=0;\\
3) \langle \varphi,\psi\rangle=\overline{\langle
\psi,\varphi\rangle};\\
4) \langle \lambda \varphi,\psi\rangle=\lambda\langle
\varphi,\psi\rangle;\\
5) \langle \varphi+\psi,\eta\rangle=\langle
\varphi,\eta\rangle+\langle \psi,\eta\rangle.$

If $\langle \cdot,\cdot\rangle:\mathcal{A}\times
\mathcal{A}\rightarrow L^0$ is an $L^0$-valued inner product then
the following formula
$$
\|\varphi\|=\sqrt{\langle \varphi,\varphi \rangle}
$$
determines an $L^0$-\textit{valued norm} on $\mathcal{A}.$ A pair
$\langle\mathcal{A},\langle \cdot,\cdot\rangle)$ is called
\textit{Kaplansky-Hilbert module} over $L^0$ or
$L^0$-\textit{Hilbert space} if $(\mathcal{A},\|\cdot\|)$ is BKS
over $L^0$ (see [8, 9]).

Let $X$ be a Kaplansky-Hilbert module over $L^0,$ and $X_0\subset
X.$ Note that $X_0$ is a $bo$-{\it closed submodule} of the
Kaplansky-Hilbert module $X$ if and only if $X_0$ is a submodule in
the usual sense, i. e. $X_0$ is a set containing all sums of the
form $bo$-$\sum\limits_{\alpha\in A}\pi_\alpha\varphi_\alpha$, where
$(\varphi_\alpha)_{\alpha\in A}$ is any bounded family in $X_0$ and
$(\pi_\alpha)_{\alpha\in A}$ is a partition of the unit in $\nabla,$
and  it is also closed with respect to the norm of the module $X.$

Let $I$ be an index set. For every $i\in I$ consider a
Kaplansky-Hilbert module $X_{i}$ over $L^0.$ Put $X_{I}=\{\varphi\in
\prod\limits_{i\in I} X_{i}: (o)$-$\sum\limits_{i\in
I}\|\varphi_{i}\|^{2}_{i} \in L^0\}.$ Considered with the pointwise
operations, $X_I$ forms a module over $L^0.$ The inner product
$\langle\cdot,\cdot\rangle:X_{I}\times X_{I}\rightarrow L^0$ is
defined as follows:
$$
\langle \varphi,\psi\rangle=(o)\mbox{-}\sum\limits_{i\in I}\langle
\varphi_{i},\psi_{i}\rangle_{i},
$$
where $\varphi,\psi\in X_{I}$ and
$\langle\cdot,\cdot\rangle_{i}:X_{i}\times X_{i}\rightarrow L^0$ is
the inner product in the corresponding $X_i$. Then
$\|\varphi\|=\sqrt{\langle{\varphi,\varphi}\rangle}$ gives an
$L^0$-valued norm on $X_I,$ and it clear that
$\|\varphi\|=((o)$-$\sum\limits_{i\in I}\langle
\varphi_{i},\varphi_{i}\rangle_{i})^{1/2}.$ Besides $X_{I}$ equipped
with this structure forms a Kaplansky-Hilbert module over $L^0.$ We
say that $X_{I}$ is the direct sum of the family $(X_{i})_{i\in I}$
and denote it by $\bigoplus\limits_{i\in I} X_{i}$.

Let $X_1,$ $X_2$ be Kaplansky-Hilbert modules over $L^0,$ and let
$a$ be an operator from $X_1$ into $X_2$. The domain of the operator
$a$ is denoted by $\mathcal{D}(a).$ The set of all pairs $(\varphi,
a\varphi),\ \varphi\in \mathcal{D}(a),$ in the direct sum $X_1\oplus
X_2,$ is called the {\it graph} of the operator $a$. The graph of
the operator $a$ is denoted by $G(a).$ Thus
$$G(a)=\{(\varphi, a\varphi):\ \varphi\in D(a)\}.$$
It is clear that two operators $a$ and $b$ coincide if and only if
$G(a)=G(b).$ The set $S\subset X_1\oplus X_2$ is the graph of an
appropriate operator if and only if the relations $(\varphi,
\psi)\in S,$ $(\varphi, \psi')\in S$ imply $\psi=\psi'$. An operator
$a:X_1\rightarrow X_2$ is $L^0$-linear if and only if $G(a)$ is a
submodule of $X_1\oplus X_2.$ An operator $a:X_1\rightarrow X_2$ is
called $bo$-{\it closed} if its graph $G(a)$ $bo$-closed in
$X_1\oplus X_2.$

If an operator $a$ is not $bo$-closed then by the definition its
graph $G(a)$ is not $bo$-closed in $X_1\oplus X_2.$ If the
$bo$-closure $\overline{G(a)}$ of the set $G(a)$ in $X_1\oplus X_2$
is the graph of some operator, then this operator is denoted by
$\widetilde{a}$ and it is called the $bo$-{\it closure} of $a$. In
this case the operator $a$ is said to be $bo$-{\it closable}
operator.

Note that $\widetilde{a}$ is the least $bo$-closed extension of the
operator $a$. The set $\overline{G(a)},$ which is the graph of the
operator $\widetilde{a}:X_1\rightarrow X_2,$ consists of elements of
the form $(\varphi, a\varphi),$ $\varphi\in \mathcal{D}(a)$ and
their $bo$-limits.

For a Kaplansky-Hilbert module $X$ over $L^0$ an $L^0$-valued
version of the Riesz theorem is also true , i. e. for every
$L^0$-bounded $L^0$-linear functional $f:X\rightarrow L^0$ there
exists a vector $\psi\in X$ such that $f(\varphi)=\langle\varphi,
\psi\rangle$ for all $\varphi\in X$ (see [9]).

Let $a:X\rightarrow Y$ be an $L^0$-linear operator. An {\it adjoint}
operator to $a$ is an operator $a^*: Y\rightarrow X$, satisfying the
condition $\langle a\varphi, \psi\rangle=\langle \varphi, a^*
\psi\rangle$ for all $\varphi\in X$ and $\psi\in Y$.

Let $\varphi,\ \psi\in X.$ We define an $L^{0}$-linear operator
$\varphi\otimes \psi$ on $X$ by the rule
$$(\varphi\otimes \psi)\eta=\langle \eta,\psi\rangle \varphi. $$

An element $\lambda\in L^{0}$ is called \textit{strictly positive}
(denoted by $\lambda\gg 0$) if $\lambda(\omega)>0 $ for almost every
$\omega\in \Omega$. If $\|\varphi\|\gg 0,\ \|\psi\|\gg 0,$ then the
operator $\varphi\otimes \psi$ is homogenous of type one. Moreover,
$\varphi\otimes \psi$ is a projection if and only if $\psi=\varphi$
and $\|\varphi\|=\mathbf{1}.$

\section*{\center 2. Derivations and automorphisms of $O^\ast$-algebras over $L^0$}

Let $X$ be a Kaplansky-Hilbert module over $L^0,$ and let
$\mathcal{D}\subset X$ be a dense domain. By $I_\mathcal{D}$ we
denote the identity map on $\mathcal{D}$.

{\bf Definition 1.} A set of $bo$-closable $L^0$-linear operators
with the domain $\mathcal{D}$ and containing $I_\mathcal{D}$ is said
to be an \textit{$O$-family} over $L^0$. In this case $\mathcal{D}$
is called the \textit{domain} of this family.

If $\mathcal{A}$ is an $O$-family over $L^0$ then the domain of this
family will be denoted by $\mathcal{D}(\mathcal{A}).$ If $a\in
\mathcal{A}$ then according to the definition we have
$\mathcal{D(A)}=\mathcal{D}(a)=\mathcal{D}.$

{\bf Definition 2.} An \textit{$O$-module} over $L^0$ is an
$O$-family $\mathcal{A}$ over $L^0$ such that $\alpha a+\beta b \in
\mathcal{A}$ for all $a,\ b \in \mathcal{A}$ and $\alpha,\ \beta \in
L^0$.

Recall that by $ab$ we denote the composition of the operators $a$
and $b$. If $a$ and $b$ are operators on $\mathcal{D}$ and
$b\mathcal{D}\subset \mathcal{D}$ then $ab$ is also an operator on
$\mathcal{D}$ defined by $ab\varphi=a(b\varphi),$ $\varphi\in
\mathcal{D}$.

{\bf Definition 3.} An \textit{$O$-algebra} over $L^0$ is an
$O$-module $\mathcal{A}$ over $L^0$ such that
$b\mathcal{D}(\mathcal{A})\subset \mathcal{D}(\mathcal{A})$ and
$ab\in \mathcal{A}$ for all $a,\ b \in \mathcal{A}$.

It is easy to see that every $O$-algebra over $L^0$ with the
operations of addition, multiplication by elements of $L^0$ and the
product defined as the composition of operators, is an algebra over
$L^0$. Note also that $I_\mathcal{D}$ is the unit of this algebra.

{\bf Definition 4.} An \textit{$O^\ast$-family} over $L^0$ on
$\mathcal{D}$ is a set $\mathcal{A}$ of $L^0$-linear operators with
the domain $\mathcal{D}$ such that $I_\mathcal{D}\in \mathcal{A},$
$\mathcal{D}\subset \mathcal{D}(a^\ast),$ and $a^+\in \mathcal{A}$
for all $a\in \mathcal{A}$, where $a^+=a^\ast|\mathcal{D}$.

Let $\mathcal{A}$ be an $O^\ast$-family over $L^0$ on $\mathcal{D}$.
Then $\mathcal{A}$ is an $O$-family over $L^0$ on $\mathcal{D}$.
Indeed, each operator $a\in\mathcal{A}$ is $bo$-closable because
$\mathcal{D}\subset\mathcal{D}(a^*)$ and $\mathcal{D}$ is dense in
$X.$

If $a\in \mathcal{A}$ then
$$
\langle a\varphi, \psi \rangle=\langle \varphi, a^+\psi \rangle \
\mbox{for all} \ \varphi,\ \psi\in \mathcal{D} \eqno (1)
$$
and hence $a=(a^+)^+.$ From the above  we obtain, in particular,
that $a\rightarrow a^+$ is a bijective map of $\mathcal{A}$ onto
itself.

{\bf Definition 5.} An \textit{$O^\ast$-module} over $L^0$ is an
$O$-module over $L^0$ which is an $O^\ast$-family over $L^0$.

If $\mathcal{A}$ is an $O^*$-module over $L^0$ on $\mathcal{D}$ then
the map $a\rightarrow a^+,$ $a\in \mathcal{A},$ is an involution on
$\mathcal{A}.$

{\bf Definition 6.} An \textit{$O^\ast$-algebra} over $L^0$ is an
$O$-algebra over $L^0$ which is an $O^\ast$-family over $L^0$.

Let $\mathcal{L^+(D)}$ denote the set of all $L^0$-linear operators
$a$ on a Kaplansky-Hilbert module $X$ over $L^0$ which satisfy
$a\mathcal{D}\subset \mathcal{D},$
$\mathcal{D}\subset\mathcal{D}(a^*)$ and $a^*\mathcal{D}\subset
\mathcal{D}.$

{\bf Theorem 1}. \textit{$\mathcal{L^+(D)}$ is the largest
$O^*$-algebra over $L^0$ with the domain $\mathcal{D}.$}

{\it Proof}. At first we check that $\mathcal{L^+(D)}$ is an
$O^*$-family over $L^0$. Let $a\in \mathcal{L^+(D)}.$ We have
$a^+\mathcal{D} =a^*\mathcal{D} \subset\mathcal{D},$
$(a^+)^*=(a^*|\mathcal{D})^*\supset a^{**}\supset a,$ and hence
$(a^+)^*\mathcal{D}=a\mathcal{D} \subset\mathcal{D},$ i. e. $a^+\in
\mathcal{L^+(D)},$ as it was asserted.

Now let us show that $\mathcal{L^+(D)}$ is an $O$-algebra over
$L^0$. Let $a, b\in \mathcal{L^+(D)}.$ It is easy to see that
$\lambda a\in \mathcal{L^+(D)}$ for all $\lambda\in L^0.$ From
$\mathcal{D}\subset\mathcal{D}(a^*)\cap\mathcal{D}(b^*)\subset\mathcal{D}((a+b)^*)$
and $(a+b)^*\mathcal{D}=(a^*+b^*)\mathcal{D}$ it follows that
$(a+b)\in\mathcal{L^+(D)}.$

We shall show that $ab\in\mathcal{L^+(D)}.$ Let $\varphi\in
\mathcal{D}$ and $\psi\in \mathcal{D}.$ According to (1) we have
$\langle ab\varphi, \psi\rangle = \langle b\varphi, a^+\psi\rangle.$
By virtue of $a^+\mathcal{D} \subset\mathcal{D},$ applying again
(1), we obtain $\langle ab\varphi, \psi\rangle = \langle \varphi,
b^+a^+\psi\rangle.$ Besides, $b^+a^+\subset (ab)^*$ and
$b^+a^+=(b^*|\mathcal{D})(a^*|\mathcal{D})=(b^*a^*)|\mathcal{D}=(ab)^*|\mathcal{D}=(ab)^+.$
These imply that $\mathcal{D}\subset\mathcal{D}((ab)^*),$ $(ab)^*
\mathcal{D}= b^+a^+ \mathcal{D}\subset \mathcal{D}.$ Thus, $ab\in
\mathcal{L^+(D)}.$

From the above it is clear that $\mathcal{L^+(D)}$ is an
$O^*$-algebra over $L^0$.

Now let $\mathcal{A}$ be an arbitrary $O^*$-algebra over $L^0$ with
the domain $\mathcal{D}$ and let $a\in \mathcal{A}.$ According to
the definition 3 we have $a\mathcal{D} \subset \mathcal{D}$ since
$\mathcal{A}$ is an $O$-algebra. The definition 4 yields that
$a^+\in\mathcal{A}$ since $\mathcal{A}$ is an $O^*$-algebra. Hence,
$a^*\mathcal{D}=a^+ \mathcal{D} \subset\mathcal{D}.$ This means that
$\mathcal{A} \subset\mathcal{L^+(D)}.$ Theorem 1 is proved.

Let $X$ be a Kaplansky-Hilbert module over $L^0,$ and let
$\mathcal{D}\subset X$ be a $bo$-dense submodule. By the symbol
$\mathcal{L(D)}$ we denote the algebra of all $L^0$-linear operators
$a:\mathcal{D}\rightarrow\mathcal{D}$. Let $\mathcal{U}$ be a
standard algebra in $\mathcal{L(D)}$. Recall that a linear operator
$\delta:\mathcal{U}\rightarrow\mathcal{L}(\mathcal{D})$ is said to
be a {\it derivation}, if $\delta(ab)=\delta(a)b+a\delta(b)$ for all
$a, b\in \mathcal{U}.$ If for a derivation
$\delta:\mathcal{U}\rightarrow\mathcal{L}(\mathcal{D})$ there exists
an element $x\in\mathcal{U}$ such that $\delta(a)=xa-ax $ for all
$a\in \mathcal{U}$ then $\delta$ is called an \textit{inner
derivation}.

Further in theorems 2 and 3 we suppose that there exists a vector
$e$ in the $bo$-dense submodule $\mathcal{D}$ of the
Kaplansky-Hilbert module $X$ over $L^0$ such that
$\|e\|=\textbf{1},$ where $\textbf{1}$ is the unit in $L^0.$

\textbf{Theorem 2.} \emph{Let
$\delta:\mathcal{U}\rightarrow\mathcal{L(D)}$ be an $L^0$-linear
derivation of a standard algebra $\mathcal{U}.$ Then there exists
$x\in\mathcal{L(D)}$ such that
$$\delta(a)=xa-ax $$
for all $a\in\mathcal{U}.$ }

{\it Proof}. At first consider the case
$\mathcal{U}=\mathcal{F}(\mathcal{D}),$ where $\mathcal{F(D)}$ is
the algebra of finite-generated operators
$a:\mathcal{D}\rightarrow\mathcal{D}$.

Fix a vector $e\in \mathcal{D}$ with $\|e\|=\textbf{1}$ and a
functional $f:\mathcal{D}\rightarrow L^0$ such that
$f(e)=\textbf{1}.$ Define a projection
$p\in\mathcal{F}(\mathcal{D})$ by
$$p(\varphi)=f(\varphi)e, \quad \varphi\in \mathcal{D}.$$ Since
$p^{2}=p$ then $\delta(p)=p\delta(p)+\delta(p)p$ and therefore
$p\delta(p)p=0.$ Put $\psi=p\delta(p)-\delta(p)p.$ Then $p\psi-\psi
p=p\delta(p)+\delta(p)p=\delta(p).$

Putting $\delta'(a)=\delta(a)-(a\psi-\psi a)$ we get $\delta'(p)=0$.
Thus, one may assume that $\delta(p)=0.$ Then we have
$$\delta(ap)=a\delta(p)+\delta(a)p=\delta(a)p. \eqno (2)$$

Consider a vector $\varphi\in \mathcal{D}$ and an operator
$a\in\mathcal{F}(\mathcal{D})$ such that $a(e)=\varphi$. Define an
operator $x: \mathcal{D}\rightarrow\mathcal{D}$ by the formula
$$x(\varphi)=\delta(a)e.$$
The operator $x$ is defined correctly. Indeed, let
$\varphi\in\mathcal{D}$ be a vector and let $a_1,\
a_2\in\mathcal{F}(\mathcal{D})$ be operators such that
$a_1(e)=a_2(e)=\varphi.$ For each $\eta\in\mathcal{D}$ we have $(a_i
p)\eta=f(\eta)a_i(e),$ $i=1,\ 2,$ i. e. $a_1p=a_2p.$ Therefore by
virtue of (2) it follows that
$\delta(a_1)(e)=(\delta(a_1)p)(e)=\delta(a_1p)(e)=\delta(a_2p)(e)=(\delta(a_2)p)(e)=\delta(a_2)(e),$
i. e. $\delta(a_1)=\delta(a_2).$

It easy to see that the operator $x$ is $L^0$-linear.

Let $\varphi\in\mathcal{D}$ and $a\in\mathcal{F}(\mathcal{D}).$
Then $(xap)\varphi=x(a(p(\varphi)))=x(f(\varphi)a(e))=f(\varphi)x(a(e))=
f(\varphi)\delta(a)(e)=\delta(a)p(\varphi)=\delta(ap)\varphi.$
Thus, $xap=\delta(a)p$ for all $a\in \mathcal{F(D)}.$
Therefore for $b\in \mathcal{F(D)}$ we have
$xabp=\delta(ab)p= a\delta(b)p+\delta(a)bp=axbp+\delta(a)bp,$ i. e.
$$\delta(a)bp=xabp-axbp. \eqno (3)$$

Now for an arbitrary $\varphi\in \mathcal{D}$ take
$b\in\mathcal{F(D)}$ such that $b(e)=\varphi.$ Then
$(bp)(e)=\varphi.$ Hence from (3) we obtain $\delta(a)=xa-ax$ for
all $a\in \mathcal{F(D)}.$

Let now $\mathcal{U}\subset \mathcal{L(D)}$ be an arbitrary standard
algebra and take $b\in \mathcal{U}.$ Then $ba\in\mathcal{F(D)}$ for
all $a\in \mathcal{F(D)}.$ Therefore
$$\delta(ba)=xba-bax. \eqno (4)$$
On the other hand according to the definition of derivation we have
$$\delta(ba)=\delta(b)a+b\delta(a)=\delta(b)a+b(xa-ax). \eqno (5)$$
From $(4)$ and $(5)$ we obtain $\delta(b)a=xba-bxa=(xb-bx)a.$

Now for an arbitrary $\varphi\in\mathcal{D}$ take
$a\in\mathcal{F(D)}$ such that $a(\varphi)=\varphi.$ Then
$\delta(b)(\varphi) = \delta(b)(a(\varphi))=(\delta(b)a)(\varphi) =
((xb-bx)a)(\varphi)= (xb-bx)(a(\varphi))=(xb-bx)(\varphi),$ i. e.
$\delta(b)(\varphi) = (xb-bx)(\varphi)$ for any
$\varphi\in\mathcal{D}.$ This means that $\delta(b)=xb-bx$ for all
$b\in \mathcal{U}.$ Theorem 2 is proved.

Replacing $\mathcal{F(D)}$ by $\mathcal{F^+(D)}:=\mathcal{F(D)} \cap
\mathcal{L^+(D)}$ and $\mathcal{L(D)}$ by $\mathcal{L^+(D)},$ we get

\textbf{Corollary 1.} \emph{Let
$\delta:\mathcal{U}\rightarrow\mathcal{L^+(D)}$ be an $L^0$-linear
derivation of the algebra $\mathcal{U}\supset\mathcal{F^+(D)},$
where $\mathcal{D}$ is a $bo$-dense submodule of a Kaplansky-Hilbert
module $X$ with a vector $e\in\mathcal{D}$ with $\|e\|=\mathbf{1}.$
Then there exists $x\in\mathcal{L^+(D)}$ such that
$$\delta(a)=xa-ax $$
for all $a\in\mathcal{U}.$ In particular each $L^0$-linear
derivation of the algebra $\mathcal{L}^+(\mathcal{D})$ over $L^0$ is
inner.}

Recall that a bijective linear operator
$\alpha:\mathcal{L(D)}\rightarrow\mathcal{L(D)}$ is called
\textit{automorphism} if $\alpha(ab)=\alpha(a)\alpha(b)$ for all
$a,\ b\in \mathcal{L(D)}.$

\textbf{Theorem 3.} \emph{Let
$\alpha:\mathcal{F(D)}\rightarrow\mathcal{F(D)}$ be an $L^0$-linear
automorphism of the algebra $\mathcal{F(D)}.$ Then there exists
$x\in\mathcal{L(D)}$ such that $x^{-1}\in\mathcal{L(D)}$ and
$$\alpha(a)=xax^{-1}$$
for all $a\in\mathcal{F(D)}.$}

{\it Proof}. Let $e\in\mathcal{D}$ be a vector with
$\|e\|=\mathbf{1}$ and let $f:\mathcal{D}\rightarrow L^0$ be an
$L^0$-linear functional such that $\|e\|=\textbf{1},\
f(e)=\textbf{1}.$ We define a projection $p\in\mathcal{F(D)}$ as
follows
$$p(\varphi)=f(\varphi)e, \quad \varphi\in \mathcal{D}.$$ Then
obviously $p(e)=e.$ Moreover the projection $\alpha(p)$ is
homogeneous of type one because $\alpha$ is an $L^0$-linear
automorphism. Now take $e_1\in \mathcal{D}$ such that
$\|e_1\|=\textbf{1},\ \alpha(p)(e_1)=e_1.$

We define an operator $x:\mathcal{D}\rightarrow \mathcal{D}$ as
follows: for any $\varphi\in \mathcal{D}$ take an operator
$a\in\mathcal{F(D)}$ such that $a(e)=\varphi$ and put
$$x(\varphi)=\alpha(a)(e_1),\quad \varphi\in \mathcal{D}.$$

Let $\varphi\in \mathcal{D}$ and take $a_1, a_2\in\mathcal{F(D)}$
such that $a_1(e)=a_2(e)=\varphi.$ For each $\psi\in\mathcal{D}$ we
have $(a_i p)(\psi)=f(\psi)a_i(e),\ i=1, 2,$ i. e. $a_1p=a_2p.$
Therefore
$\alpha(a_1)(e_1)=\alpha(a_1)\alpha(p)(e_1)=\alpha(a_1p)(e_1)=\alpha(a_2p)(e_1)=\alpha(a_2)\alpha(p)(e_1)
=\alpha(a_2)(e_1).$ This means that $x$ is defined correctly.

Obviously $x$ is $L^0$-linear.

Now we shall show that $x$ is a bijection. Let $\varphi_1,\
\varphi_2\in \mathcal{D}$ such that $\varphi_1\neq \varphi_2.$
Choose $a_1,\ a_2\in\mathcal{F(D)}$ such that $a_i(e)=\varphi_i,\
i=1,\ 2.$ Then $a_1\neq a_2,$ and hence $a_1p\neq a_2p.$ Since
$a_ip,$ $i=1, 2,$ are one-generated operators and $\alpha$ is an
automorphism then $\alpha(a_1)(e_1)= \alpha(a_1)\alpha(p)(e_1) =
\alpha(a_1 p)(e_1) \neq \alpha(a_2 p)(e_1)
=\alpha(a_2)\alpha(p)(e_1)= \alpha(a_2)(e_1).$ Hence,
$x(\varphi_1)\neq x(\varphi_2).$ Now take $\psi\in\mathcal{D},$ and
$a\in\mathcal{F(D)}$ such that $a(e_1)=\psi.$ Put
$b=\alpha^{-1}(a).$ Then for $\varphi=b(e)$ one has
$x(\varphi)=\alpha(b)(e_1)=\alpha(\alpha^{-1}(a))(e_1)=a(e_1)=\psi,$
i. e. $x(\varphi)=\psi.$

Let $\varphi\in \mathcal{D}$ and $a\in\mathcal{F(D)}.$ Take
$b\in\mathcal{F(D)}$ such that $b(e)=\varphi.$ Then
$(xa)(\varphi)=x(a(\varphi))= x(ab(e)) =
\alpha(ab)(e_1)=\alpha(a)\alpha(b)(e_1)=\alpha(a)x(\varphi).$ Thus,
$xa=\alpha(a)x,$ i. e. $\alpha(a)=xax^{-1}$ for all
$a\in\mathcal{F(D)}.$ Theorem 3 is proved.

\textbf{Corollary 2.} \emph{For each $L^0$-linear automorphism of a
standard algebra $\mathcal{U}$ there exists $x\in\mathcal{L(D)}$
such that $x^{-1}\in\mathcal{L(D)}$ and
$$\alpha(a)=xax^{-1}$$
for all $a\in\mathcal{U}.$ In particular, each $L^0$-linear
automorphism of the algebra $\mathcal{L(D)}$ is spatial.}

Let $\mathcal{D}_1$, resp. $\mathcal{D}_2$ be $(bo)$-dense
submodules in the Kaplansky-Hilbert modules $X_1$, resp. $X_2$ over
$L^0$, and let $\mathcal{A}_1$ and $\mathcal{A}_2$ be
$*$-subalgebras respectively in the $O^*$-algebras
$\mathcal{L}^+(\mathcal{D}_1)$ and $\mathcal{L}^+(\mathcal{D}_2)$
over $L^0$.

\textbf{Definition 7.} An $L^{0}$-linear $*$-isomorphism
$\pi:\mathcal{A}_1{\longrightarrow} \mathcal{A}_2$ is said to be
\textit{spatial} if there exists an isometry
$U:X_1\stackrel{\mbox{on}}{\longrightarrow} X_2$ such that

(i) $U\mathcal{D}_1=\mathcal{D}_2,$

(ii)$\pi(a)\varphi=UaU^{-1}\varphi \ \mbox{for all} \
\varphi\in\mathcal{D}_2,\ a\in\mathcal{A}_1.$\\
Then we say that $\pi$ is \textit{implemented} by the operator $U$.

An $L^0$-linear $\ast$-automorphism of an algebra $\mathcal{A}$  is
called \textit{inner}, if it is spatial and it may be implemented by
a unitary operator $U$ on a Kaplansky-Hilbert module $X$ over $L^0$
such that $U|\mathcal{D}\in\mathcal{A},$ where $\mathcal{D}$ is a
$(bo)$-dense submodule of $X$.

Let $A$ be a module over $L^0$ and a $*$-algebra over $L^0$. The set
of all projections in $A$ is denoted by $I_{sa}(A).$ If $p_1,p_2\in
I_{sa}(A)$ then we write $p_1 \leq\ p_2$ if and only if $p_1
p_2=p_1$. The relation $\leq$ is a reflexive, antisymmetric and
transitive relation in $I_{sa}(A)$. If there exists an
$L^{0}$-linear $*$-isomorphism $\pi$ from the algebra $A$ onto a
$*$-subalgebra of $\mathcal{L}^{+}(\mathcal{D})$, and if $p$ is a
projection in $A$ then $\widetilde{\pi(p)}$ is also a projection in
$\mathcal{L}^{+}(\mathcal{D})$, i. e.
$\widetilde{\pi(p)}\in\mathcal{L}^{+}(\mathcal{D})$ and
$\widetilde{\pi(p)}^{2}=\widetilde{\pi(p)}.$ Obviously the relation
$p_1\leq\ p_2$ is equivalent to the usual relation
$\widetilde{\pi(p_1)}\leq\widetilde{\pi(p_2)}$ between the
projections $\widetilde{\pi(p_1)}$ and $\widetilde{\pi(p_2)}$. Let
$H_{1}(A)$ denote the set of all homogeneous of type one projections
of the algebra $A.$ For $p_1,p_2\in H_{1}(A)$ we shall write
$p_1\approx p_2,$ if $p_1 Ap_2\neq \{0\}.$ Further on, the elements
of the set $H_{1}(A)$ will be called projections of rank one.

Let $\mathcal{D}_{i}$ be a $(bo)$-dense submodule of a
Kaplansky-Hilbert module $X_i$ over $L^0$ such that there exists
$\varphi_{i}\in\mathcal{D}_{i},$ $\|\varphi_{i}\|=\mathbf{1},$ $i\in
I.$ By $\mathcal{D}_I$ we denote a $(bo)$-dense submodule of the
Kaplansky-Hilbert module $X_I$ over $L^0$, consisting of all vectors
$(\varphi_i):=(\varphi_i)_{i\in I},$ which have only finitely many
nonzero coordinates $\varphi_i\in\mathcal{D}_i$.

Note that every element $(a_i):=(a_i)_{i\in I}$ of the product
$\prod\limits_{i\in I}\mathcal{L}^{+}(\mathcal{D}_{i})$ is an
operator on $\mathcal{D}_I$ which acts according to the formula:
$$
(a_i)(\varphi_i)=(a_i\varphi_i), \qquad (\varphi_i)\in\mathcal{D}_I.
$$
The set of all such operators forms an $O^*$-algebra with the domain
$\mathcal{D}_I.$ This algebra is denoted by
$\mathcal{L}^{+}(\mathcal{D}_{i}:i\in I)$

\textbf{Lemma 1.} \textit{Let $\mathcal{A}$ be a $*$-subalgebra of
the algebra $\mathcal{L}^{+}(\mathcal{D}_I)$ over $L^0$ and let
$M(\mathcal{A})$ be the set of all projections $p\in
H_{1}(\mathcal{A}),$ for which the generators of the images
$p\mathcal{D}_{I}$ have a unique nonzero coordinate. Then:}

\textit{(i) The set $M(\mathcal{L}^{+}(\mathcal{D}_{i}:i\in I))$
consists of the projections of the form
$\varphi_{i}\otimes\varphi_{i}$, where
$\varphi_{i}\in\mathcal{D}_{i},$ $\|\varphi_{i}\|=\mathbf{1},$ $i\in
I.$ If $\varphi_{i}\otimes\varphi_{i}$ and $\psi_{j}\otimes\psi_{j}$
are two such operators then $\varphi_{i}\otimes\varphi_{i}\approx
\psi_{j}\otimes\psi_{j}$ if and only if $i=j.$}

\textit{(ii) $M(\mathcal{A})=M(\mathcal{L}^{+}(\mathcal{D}_{i}:i\in
I))$ if and only if
$\mathcal{A}\subseteq\mathcal{L}^{+}(\mathcal{D}_{i}:i\in I)$ and
$\mathcal{F}(\mathcal{D}_{i})\subseteq\mathcal{A}$ for all $i\in
I$.}

\textit{(iii) If
$M(\mathcal{A})=M(\mathcal{L}^{+}(\mathcal{D}_{i}:i\in I))$ then
\textit{on the set} $M(\mathcal{A})=
M(\mathcal{L}^{+}(\mathcal{D}_{i}:i\in I))$ the relation "$\approx$"
corresponding to the $*$-algebra $\mathcal{A}$ coincides with the
relation "$\approx$" corresponding to the $*$-algebra
$\mathcal{L}^{+}(\mathcal{D}_{i}:i\in I)$.}

\textit{(iv) The set $H_{1}(\mathcal{L}^{+}(\mathcal{D}_{i}:i\in
I))$ of projections of rank one consists of all $L^{0}$-linear
projections of the form $(o)$-$\sum\limits_{i\in
I}\pi_i(\varphi_{i}\otimes\varphi_{i})$, where
$\varphi_{i}\in\mathcal{D}_{i},$ $\|\varphi_{i}\|=\mathbf{1},$ and
$(\pi_i)_{i\in I}$ is a partition of the unit in $\nabla.$}

\textit{Proof.} (i) From the definition it follows that the
operators of the form $\varphi_{i}\otimes\varphi_{i},$
$\varphi_{i}\in\mathcal{D}_{i}$, $\|\varphi_{i}\|=\mathbf{1}$, $i\in
I,$ are projections of rank one.

Let $\varphi_{i}\otimes \varphi_{i}, \psi_{j}\otimes \psi_{j}\in
M(\mathcal{L}^{+}(\mathcal{D}_{i}:i\in I)).$ If $i\neq j$ then
$\varphi_{i}\otimes\varphi_{i}M(\mathcal{L}^{+}(\mathcal{D}_{i}:i\in
I)) \psi_{j}\otimes \psi_{j}=\{0\} $. This implies that
$\varphi_{i}\otimes\varphi_{i}\approx \psi_{j}\otimes\psi_{j}$ if
and only if $i=j.$

ii) Suppose that
$M(\mathcal{A})=M(\mathcal{L}^{+}(\mathcal{D}_{i}:i\in I)$. At first
we shall prove that $\mathcal{A}\subseteq \mathcal{L}^{+}
(\mathcal{D}_{i}:i\in I)$. Fix $i\in I.$ If we prove that
$a\varphi\in \mathcal{D}_{i}$ for some $\varphi\in \mathcal{D}_{i}$
then by virtue of the linearity of the operator $a\in \mathcal{A}$
we have $a\psi\in \mathcal{D}_{i}$ for any $\psi\in
\mathcal{D}_{i}.$ Therefore without loss of generality we may
suppose that $\|\varphi\|=\mathbf{1}$ and $a\varphi\neq 0$. Then
$\varphi\otimes\varphi\in \mathcal{A}$ and hence
$a+\varphi\otimes\varphi\in \mathcal{A}.$ Apply the operator
$a+\varphi\otimes\varphi$ to the element $\varphi\in
\mathcal{D}_{i}:\ \ \ (a+\varphi\otimes\varphi)
(\varphi)=a(\varphi)+\varphi.$ This implies  that
$a\varphi\in\mathcal{D}_{i}.$

Now let us show that
$\mathcal{F}(\mathcal{D}_{i})\subset\mathcal{A}.$ For this it is
enough to prove that $\varphi\otimes\psi\in \mathcal{A}$ for all
unit elements $\varphi,\ \psi \in\mathcal{D}_{i}$ since each
finite-generated operator from $\mathcal{F}(\mathcal{D}_{i})$ may be
represented as a linear combination of operators of rank one. Let
$\varphi,\psi \in \mathcal{D}_{i}$ and
$\|\varphi\|=\|\psi\|=\mathbf{1}.$ By virtue of (i) we have
$\varphi\otimes\varphi, \psi\otimes\psi \in
M(\mathcal{A})=M(\mathcal{L}^{+}(\mathcal{D}_{i}:i\in I)).$ From
this it follows that the operators $\varphi\otimes\varphi$,
$\psi\otimes\psi$ belong to $\mathcal{A}$ and hence
$$(\psi\otimes\psi)(\varphi\otimes\varphi)=\langle
\psi,\varphi\rangle (\varphi\otimes\psi)\in \mathcal{A},$$ i. e.
$(\varphi\otimes\psi)\in \mathcal{A}.$

The inverse statement is obvious.

(iii) If $M(\mathcal{A})=M(\mathcal{L}^{+}(\mathcal{D}_{i}:i\in I))$
then from (ii) it follows that $\mathcal{A}\subset \mathcal{L}^{+}
(\mathcal{D}_{i}:i\in I).$ Therefore according to (i) it is
sufficient to show that
$\psi\otimes\psi\mathcal{A}\varphi\otimes\varphi\neq\{0\}$ for all
unit elements $\varphi,\psi\in\mathcal{D}_{i}.$ Consider
$\xi\otimes\xi\in \mathcal{A},$ where the vector $\xi\in
\mathcal{D}_{i}$ is defined by the formula
$$
\xi= \left\{
\begin{array}{ll}
\frac{1}{\sqrt{2}}(\varphi+\psi) &
\mbox{if}\ \langle\varphi,\ \psi\rangle=0,\\
\ \varphi &\mbox{in other cases.}
\end{array}
\right.
$$
Then we have
$$
(\psi\otimes\psi)(\xi\otimes\xi)(\varphi\otimes\varphi)=\langle\varphi,\xi\rangle
\langle\psi,\xi\rangle \varphi\otimes\psi\neq 0,
$$
i. e. $\psi\otimes\psi\mathcal{A}\varphi\otimes\varphi\neq\{0\}$.

(iv) Let $a=(a_{i})_{i \in I}\in \mathcal{L}^{+}(\mathcal{D}_{i}:i
\in I)$ be a projection of rank one. Then $a_{i}$ is a projection in
$\mathcal{L}^{+}(\mathcal{D}_{i})$ for all $i \in I.$ Since $a$ is a
projection of rank one there exist a partition $(\pi_i)_{i\in I}$ of
the unit in $\nabla$ and a vector $\varphi_{i}\in \mathcal{D}_{i},\
\ \|\varphi_{i}\|=\mathbf{1},$ such that
$a_{i}=\pi_i(\varphi_{i}\otimes\varphi_{i}).$ From this we have
$a=(o)$-$\sum\limits_{i\in I}\pi_i(\varphi_{i}\otimes\varphi_{i}).$
Lemma 1 is proved.

\textbf{Theorem 4.} \textit{Let $\mathcal{D}_i$ and
$\mathcal{D}_{j}$ be $(bo)$-dense submodules of Kaplansky-Hilbert
modules $X_i$ ($i\in I$) and $X_j$ ($j\in J$) over $L^0,$
respectively, such that for each $i\in I$ and $j\in J$ there exist
$e_i\in\mathcal{D}_i$ and $f_j\in\mathcal{D}_{j}$ with
$\|e_i\|=\mathbf{1}$ and $\|f_j\|=\mathbf{1}$. Let $\mathcal{A}$ and
$\mathcal{B}$ be $*$-subalgebras of the algebras
$\mathcal{L}^{+}(\mathcal{D}_I)$ and
$\mathcal{L}^{+}(\mathcal{D}_J)$ over $L^0$, respectively,
satisfying the following conditions
$$M(\mathcal{A})=M(\mathcal{L}^{+}(\mathcal{D}_{i}:i\in I)),$$
$$M(\mathcal{B})=M(\mathcal{L}^{+}(\mathcal{D}_{j}:j\in J)).$$
Suppose that there exists an $L^0$-linear $*$-isomorphism $\pi,$
mapping $\mathcal{A}$ onto $\mathcal{B}$. Then $\pi$ is a spatial
$L^0$-linear $*$-isomorphism. Moreover, there exist a partition
$(\pi_{\alpha})$ of the unit in $\nabla,$ bijective maps
$\chi_{\alpha}:I\rightarrow J$ and surjective isometries
$U_{\alpha}:X_{I}\rightarrow X_{J}$ such that
$U=\sum\limits_{\alpha}\pi_{\alpha}U_{\alpha}$ implements  $\pi$ and
$U_{\alpha}(\pi_{\alpha}\mathcal{D}_{i})=\pi_{\alpha}\mathcal{D}_{\chi_{_{\alpha}}(i)}$
for all $i\in I$.}

\textit{Proof.} Since $\pi$ is a $*$-isomorphism, it preserves the
relation $\approx$ and $\pi(M(\mathcal{A}))\subset
H_{1}(\mathcal{B}).$ Hence
$$
\pi(M(\mathcal{L}^{+}(\mathcal{D}_{i}:i\in I)))\subset
H_{1}(\mathcal{L}^{+}(\mathcal{D}_{j}:j\in J)). \eqno(6)
$$

From (6) we have $\pi(\varphi_{i}\otimes\varphi_{i})\in
H_{1}(\mathcal{L}^{+} (\mathcal{D}_{j}:j\in J)).$ By virtue of lemma
1, $\pi(\varphi_{i}\otimes\varphi_{i})$ has the form
$(o)$-$\sum\limits_{j\in J}\pi_{ij}(\psi_{ij}\otimes \psi_{ij}),$
where $(\pi_{ij})_{j\in J}$ is a partition of the unit in $\nabla$
such that $(\pi_{ij})_{i\in I}$ is also a partition of the unit in
$\nabla$.

Since $\pi$ is a $*$-isomorphism the cardinalities of the sets $I$
and $J$ are equal. Let $S(I,J)$ be the set of all bijections from
$I$ onto $J.$ For each $\alpha\in S(I,J)$ put
$\chi_{\alpha}(i)=\alpha(i)$ and
$\pi_{\alpha}=\bigwedge\limits_{i\in I}\pi_{i\chi_{_{\alpha}}(i)}.$
Then $\pi_{\alpha}\pi_{\alpha^{'}}=0$ at $\alpha\neq\alpha^{'}$ and
$\bigvee\limits_{\alpha}\pi_{\alpha}=\mathbf{1}.$ Indeed, if
$\alpha\neq\alpha^{'}$ then there exists $i_{0}\in I$ such that
$\alpha(i_{0})\neq\alpha^{'}(i_{o}).$ Then
$\pi_{i_{0}\chi_{_{\alpha}}(i_{0})}\pi_{i_{0}\chi_{_{\alpha^{'}}}(i_{0})}=0.$
From this it follows that $\pi_{\alpha}\pi_{\alpha^{'}}=0$ at
$\alpha\neq\alpha^{'}.$ Further,
$\bigvee\limits_{\alpha}\pi_{\alpha}=\bigvee\limits_{\alpha}(\bigwedge\limits_{i\in
I}\pi_{i\chi_{_{\alpha}}(i)})=\bigwedge\limits_{i\in
I}(\bigvee\limits_{\alpha\pi_{i\chi_{_{\alpha}}(i)}})=\mathbf{1}.$

Suppose that $\varphi_{i}\in\mathcal{D}_{i}$,
$\psi_{\chi_{_{\alpha}}(i)}\in\mathcal{D}_{j}$, are unit elements
such that
$\pi(\pi_{\alpha}(\varphi_{i}\otimes\varphi_{i}))=\pi_{\alpha}(\psi_{\chi_{_{\alpha}}(i)}
\otimes\psi_{\chi_{_{\alpha}}(i)}).$ We shall prove that
$$
\|\pi_{\alpha}x\varphi_{i}\|=\|\pi_{\alpha}\pi(x)\psi_{\chi_{_{\alpha}}(i)}\|
\eqno (7)
$$
for any $x\in\mathcal{A}.$ From the lemma 1 it follows that
$x\varphi_{i}\in\mathcal{D}_{i}$ and hence
$\pi_{\alpha}x\varphi_{i}\in\mathcal{D}_{i},$
$\pi_{\alpha}(x\varphi_{i}\otimes x\varphi_{i})\in \mathcal{A}$. One
has
$$\pi(\pi_{\alpha}(x\varphi\otimes x\varphi))=
\pi_{\alpha}\pi(x(\varphi\otimes\varphi)x^{+})=$$
$$\pi_{\alpha}\pi(x)\pi(\varphi\otimes\varphi)\pi(x)^{+}=
\pi(x)\pi(\pi_{\alpha}(\varphi\otimes\varphi))\pi(x)^{+}=$$
$$=\pi_{\alpha}(\pi(x)\psi_{\chi_{_{\alpha}}(i)}\otimes\pi(x)\psi_{\chi_{_{\alpha}}(i)}).
\eqno (8)$$ If $\pi(x)\psi=0$ then (7) is true. If
$\pi(x)\psi_{\chi_{_{\alpha}}(i)}\neq 0$ then
$$
(\pi(\pi_{\alpha}(x\varphi_{i}\otimes
x\varphi_{i}))^{2}=\pi(\pi_{\alpha}(x\varphi\otimes\
x\varphi)^{2})=\pi_{\alpha}\pi\|x\varphi_{i}\|^{2}(x\varphi_{i}\otimes\
x\varphi_{i})=$$
$$=\|\pi_{\alpha}x\varphi\|^{2}(\pi(x)\psi_{\chi_{_{\alpha}}(i)}\otimes\pi(x)\psi_{\chi_{_{\alpha}}(i)}).
\eqno (9) $$ On the other hand according to (8) we have
$$
(\pi(x\varphi_{i}\otimes x\varphi_{i}))^{2}=
\|\pi_{\alpha}\pi(x)\psi_{\chi_{_{\alpha}}(i)}\|^{2}(\pi(x)
\psi_{\chi_{_{\alpha}}(i)}\otimes \pi(x)\psi_{\chi_{_{\alpha}}(i)}).
\eqno (10)
$$
From the equalities (9) and (10) we obtain (7). If $i\in I$ then
from (7) it follows that the equality $$U_{\alpha
i}(\pi_{\alpha}x\varphi_{i})=\pi_{\alpha}\pi(x)\psi_{\chi_{_{\alpha}}(i)},
\ \ x\in\mathcal{A},$$ defines a unique norm preserving
$L^{0}$-linear surjective map $U_{\alpha i}:
\pi_{\alpha}\mathcal{A}\varphi_{i} \rightarrow
\pi_{\alpha}\pi(\mathcal{A})\psi_{\chi_{_{\alpha}}(i)} \equiv
\pi_{\alpha}\mathcal{B}\psi_{\chi_{_{\alpha}}(i)}$. By virtue of
lemma 1 the inclusions
$\mathcal{F}(\mathcal{D}_{i})\subseteq\mathcal{A}|\mathcal{D}_{i}
\subseteq\mathcal{L}^{+}(\mathcal{D}_{i})$ are true. From this it
follows that
$\pi_{\alpha}\mathcal{A}\varphi_{i}=\pi_{\alpha}\mathcal{D}_{i}.$
Similarly, $\pi_{\alpha}\mathcal{B}\psi_{\chi_{_{\alpha}}(i)}=
\pi_{\alpha}\mathcal{D}_{\chi_{_{\alpha}}(i)}.$ Thus,
$U_{\alpha}(\pi_{\alpha}\mathcal{D}_{i})=\pi_{\alpha}\mathcal{D}_{\chi_{_{\alpha}}(i)},$
where $U_{\alpha}=\bigoplus\limits_{i\in I}U_{\alpha i}.$ Since the
index $i\in I$ is arbitrary it follows that
$U_{\alpha}(\pi_{\alpha}\mathcal{D}_{I})=\pi_{\alpha}\mathcal{D}_{J}.$
Put $U=(o)$-$\sum\limits_{\alpha}\pi_{\alpha}U_{\alpha}.$ It is
clear that $U$ is a surjective isometry from $X_{I}$ onto $X_{J}$.
For $a\in \mathcal{A}$ one has
$$\pi(a)(\pi_{\alpha}\pi(x)\psi_{\chi_{_{\alpha}}(i)})=\pi_{\alpha}\pi(ax)\psi_{\chi_{_{\alpha}}(i)}=
U_{\alpha}\pi_{\alpha}ax\varphi_{i}=\ \
U_{\alpha}aU^{-1}_{\alpha}(\pi_{\alpha}\pi(x)\psi_{\chi_{_{\alpha}}(i)}),$$
i. e.
$\pi(a)(\pi_{\alpha}\pi(x)\psi_{\chi_{_{\alpha}}(i)})=U_{\alpha}aU^{-1}_{\alpha}
(\pi_{\alpha}\pi(x)\psi_{\chi_{_{\alpha}}(i)})$ for all
$x\in\mathcal{A},$ $\pi_{\alpha}$ and $i\in I$. The latter equality
implies $\pi(a)\psi=UaU^{-1}\psi$ for all $\psi\in\mathcal{D}_{J}$
and $a\in\mathcal{A}.$ Thus, $\pi$ is spatial and it is implemented
by $U.$ Theorem 4 is proved.

\textbf{Corollary 3.} \textit{Let $\mathcal{D}_{i}$ and
$\mathcal{D}_{j}$ be $(bo)$-dense submodules of Kaplansky-Hilbert
modules $X_i$ ($i\in I$) and $X_j$ ($j\in J$) over $L^{0},$
respectively, such that for each $i\in I$ and $j\in J$ there exist
$e_i\in\mathcal{D}_i$ and $f_j\in\mathcal{D}_{j}$ with
$\|e_i\|=\mathbf{1}$ and $\|f_i\|=\mathbf{1}$. If $\pi$ is an
$L^{0}$-linear $*$-isomorphism from
$\mathcal{L}^{+}(\mathcal{D}_{i}:i\in I)$ onto a $*$-subalgebra of
$\mathcal{L}^{+}(\mathcal{D}_J)$ such that
$M(\pi(\mathcal{L}^{+}(\mathcal{D}_{i}:i\in I)))=
M(\mathcal{L}^{+}(\mathcal{D}_{j}:j\in J))$ then $\pi$ is a spatial
$L^{0}$-linear $*$-isomorphism from
$\mathcal{L}^{+}(\mathcal{D}_{i}:i\in I)$ onto
$\mathcal{L}^{+}(\mathcal{D}_{j}:j\in J)$.}

\textit{Proof.} Assume that
$\mathcal{A}=\mathcal{L}^{+}(\mathcal{D}_{i}:i\in {I})$ and
$\mathcal{B}=\pi(\mathcal{A}).$ Then according to theorem 4 $\pi$ is
spatial. By the properties of the isometry $U$ listed in theorem 4
the map $a\mapsto UaU^{-1}$ is a surjection from
$\mathcal{L}^{+}(\mathcal{D}_{i}:i\in{I})$ onto
$\mathcal{L}^{+}(\mathcal{D}_{j}:j\in{J})$. The equality
$\pi(a)=UaU^{-1}$ implies that
$\pi(\mathcal{A})=\mathcal{L}^{+}(\mathcal{D}_{j}:j\in{J}).$
Corollary 3 is proved.

\textbf{Corollary 4.} \textit{Let $\mathcal{D}$ be a $bo$-dense
submodule of a Kaplansky-Hilbert module $X$ over $L^{0}$ such that
there exists $e\in \mathcal{D}$ with $\|e\|=\mathbf{1}$. Then each
$L^{0}$-linear $*$-automorphism of the $O^*$-algebra
$\mathcal{L}^{+}(\mathcal{D})$ is inner.}

\textit{Proof.} Put
$\mathcal{A}=\mathcal{B}=\mathcal{L}^{+}(\mathcal{D})$ and apply
Theorem 4. Then every $L^{0}$-linear $*$-automorphism $\pi$ of the
algebra $\mathcal{L}^{+}(\mathcal{D})$ is spatial. If $\pi$ is
implemented by some $U$ then $U\mathcal{D}=\mathcal{D}$ and
$U^{*}\mathcal{D}=\mathcal{D}$. So $U|\mathcal{D}\in
\mathcal{L}^{+}(\mathcal{D})$ and therefore by definition 7 $\pi$ is
inner. Corollary 4 is proved.

\vspace{1cm}

\textbf{Acknowledgments.} \emph{The second and third named authors
would like to acknowledge the hospitality of the $\,$ "Institut
f\"{u}r Angewandte Mathematik",$\,$ Universit\"{a}t Bonn (Germany).
This work is supported in part by the DFG 436 USB 113/10/0-1 project
(Germany) and the Fundamental Research  Foundation of the Uzbekistan
Academy of Sciences.}

\end{document}